\documentclass[11pt, numeric]{article}
\usepackage[margin=2cm]{geometry}
\usepackage[utf8]{inputenc}
\usepackage{helvet}
\usepackage{courier}
\usepackage{type1cm}
\usepackage{makeidx}
\usepackage{graphicx}
\usepackage{xcolor}
\usepackage{multicol}
\usepackage{multirow}
\usepackage[bottom]{footmisc}
\usepackage{booktabs}
\usepackage{authblk}
\usepackage{amsmath}
\allowdisplaybreaks[4]
\usepackage{url}
\usepackage{booktabs}
\usepackage[T1]{fontenc}
\usepackage[inline,shortlabels]{enumitem}
\usepackage{eqnarray}
\usepackage{amssymb}

\usepackage{amsthm}
\usepackage{tabularx}
\graphicspath{{images/}}
\usepackage{caption}
\usepackage{subcaption}
\usepackage[english]{babel}
\usepackage[misc]{ifsym}
\usepackage{orcidlink}
\usepackage{comment}
\usepackage{rotating}
\usepackage{todonotes}

\tikzstyle{process} = [rectangle, minimum width=3cm, minimum height=1cm,anchor=west,draw=black, fill=white]
\tikzstyle{arrow} = [thick,->,>=stealth]

\tikzstyle{process} = [rectangle, minimum width=3cm, minimum height=1cm,anchor=west,draw=black, fill=white]
\tikzstyle{arrow} = [thick,->,>=stealth]
\usetikzlibrary{calc, arrows.meta, intersections, patterns, positioning, shapes.misc, fadings, through,decorations.pathreplacing, positioning}

\tikzstyle{descript} = [text = black,align=center, minimum height=1.8cm, align=center, outer sep=0pt,font = \footnotesize]
\tikzstyle{activity} =[align=center,outer sep=1pt]

\usepackage[fitting]{tcolorbox}
\newtcboxfit{\mybox}{height=6.5cm,boxsep=0mm,top=2mm,bottom=2mm,left=2mm,right=2mm,
 nobeforeafter,width=\linewidth}

\usepackage{natbib}
\setcitestyle{numbers,square, sort&compress, semicolon}

\begin{document}
\thispagestyle{empty}

\title{\textbf{A disruption/restoration-based approach for elective surgical scheduling in a Children's Hospital}}

\author[a,b,*]{Martina Doneda \orcidlink{0000-0002-1660-3100}}
\author[c,d]{Gloria Pelizzo \orcidlink{0000-0002-5253-8828}}
\author[c]{Sara Costanzo \orcidlink{0000-0002-4367-7607}}
\author[a]{Giuliana Carello \orcidlink{0000-0001-5163-9865}}

\affil[a]{\footnotesize Politecnico di Milano, Department of Electronics, Information and Bioengineering, Milan, Italy}
\affil[b]{National Research Council of Italy, Institute for Applied Mathematics and Information Technologies, Milan, Italy }

\affil[c]{Pediatric Surgery Department, Buzzi Children's Hospital, Milan, Italy}

\affil[d]{Department of Biomedical and Clinical Science, University of Milan, Milan, Italy}
\affil[*]{Corresponding author: martina.doneda@polimi.it}

\maketitle

\noindent Currently under review for publication in the proceedings of the 5\textsuperscript{th} International Conference on Health Care Systems Engineering (HCSE 2023).
\thispagestyle{empty}

\setlength{\parindent}{0em}

\begin{abstract}
\noindent We consider the problem of scheduling elective surgeries in a Children’s Hospital, where disruptions due to emergencies and no-shows may arise. We account for two features that occur in many pediatric settings: \textit{i}) that it is not uncommon for pediatric patients to fall ill on the very day of their operation and, consequentially, to be unable to undergo surgery and \textit{ii}) that operating rooms normally reserved for elective surgeries can be used to treat emergency cases. Elective surgeries are scheduled taking into account the time spent on the waiting list and the patient's priority, which considers the severity of their condition and their surgical deadline, generating a nominal schedule. This schedule is optimized in conjunction with a series of back-up schedules: in fact, back-up schedules shall be available in advance so as to guarantee that the operating rooms activity immediately recovers in case of a disruption.\\
We propose an Integer Linear Programming-based approach for the problem. As there is no consolidated data on the features of both emergencies and no show, we enumerate a representative subset of the possible emergency and no-show scenarios and for each of them a back-up plan is designed. The approach reschedules patients in a way that minimizes disruption with respect to the nominal schedule and applies an as-soon-as-possible policy in case of emergencies to ensure that all patients receive timely care. The approach shows to be effective in managing disruptions, ensuring that the waiting list is managed properly, with a balanced mix of urgent and less urgent patients.\\
Therefore, the approach provides an effective solution for scheduling patients in a pediatric hospital, taking into account the unique features of such facilities.
\end{abstract}


\section{Introduction and literature positioning}
The management of operating theaters is one of the most studied topics in Operations Research applied to healthcare. In this work we deal with the problem of scheduling elective surgeries in a Children's Hospital considering disruptions. We consider the case of V. Buzzi Children's Hospital in Milan, a major pediatric center in Northern Italy. Our interest is motivated by the fact that scheduling pediatric surgeries presents two main features that constitute specific challenges that are not normally found in adult surgery:
 \begin{enumerate}
 \item oftentimes, pediatric hospitals do not have an operating room (OR) entirely devoted to emergencies;
 \item there is a much higher rate of no-shows than in adults, especially because of sudden fevers or other acute illnesses that may make it necessary to postpone the surgery \cite{boudreau2011}. 
 \end{enumerate}
Given a waiting list and a certain surgical capacity, the research question that guided us was how to best schedule patients considering possible disruptions, and, in case something actually disrupts the plan, how to best manage the rescheduling process. In this rather dynamic environment, scheduling elective surgeries is a complex challenge. This work explores a novel approach to scheduling that takes into account these disruptions, specifically tailored to the unique demands of children's healthcare facilities.\\

There is immense literature on OR planning \citep{cardoen2010operating, rothstein2018operating, zhu2019operating}. The interested reader is referred to the cited sources for a more comprehensive review of the problem of OR management. 
Commonly, this literature is divided in three tiers, corresponding to three different decision levels:
\begin{enumerate}
 \item at the strategic level, the planning scope is long, oftentimes structural, for instance dealing with how many ORs are needed in the hospital \citep{blake2002goal}. In this work, this level is not studied, and OR capacity is assumed to be given and invariable;
 \item at the tactical level, decisions deal with medium-term problems such as how to allocate ORs to different specialties \citep{dexter2005tactical}. In this work, it is assumed that the capacity of ORs is shared among surgical teams, and we assume an open scheduling policy \citep{marques2019}, with each OR and day treated as a unicum. Within this strategy, the compatibility of a surgery with an OR is defined for each day;
 \item the operational level deals with short-term decisions, one above all when to schedule each surgery \citep{erdogan2011surgery}. This is the level in which our work positions itself.
\end{enumerate}

The operational level can further be divided in two sequential processes: \textit{advance scheduling} and \textit{allocation scheduling}, the former being the allocation of a surgery to a certain day, and the latter one determining the operating room and the starting time of the procedure. In this work, we address both.\\

Pediatric scheduling presents its own unique features that make it stand apart from its adult counterpart. Planners have longtime been concerned about same-day cancellations, rescheduling and last-minute change of plans, as it causes extensive distress both to patients and their families \citep{tait1997, pratap2015}. In the considered problem, these disruptions can be due to the unforeseeable possibility of an emergency occurring during the surgical working hours or to one of the patients scheduled for surgery being unable to undergo it. Within this context, we proactively take into account unforeseen events that may disrupt the original planning \citep{denton2007optimization} by approximating them with scenarios. We use a disruption/restoration approach, developed in the field of telecommunications and transportation planning \citep{lee2007restoration}. \\
As mentioned above, the disruptions we take into account are of two kinds: emergencies, modelled as external inputs that may impact the planned schedule by delaying subsequent operations and occupying OR capacity, and no-shows, defined as the last-minute cancellation of a planned operations.\\

Accounting for emergencies is, by definition, accounting for events that {cannot} be predicted. Many strategies have been proposed to tackle them. In some works, a OR is devoted exclusively to manage any emergency that may arise \cite{aringhieri2017}, but this configuration can become undesirable when there are not many ORs and/or the number of emergencies is relatively small. Other works suggest sharing OR capacity among elective surgeries and emergencies. This is mainly done in two ways: \textit{i}) reserving some OR capacity to accommodate for emergencies \cite{rachuba2017} and \textit{ii}) not explicitly reserving any capacity for emergencies, but proactively scheduling elective procedures in a manner that facilitates emergency management.
There are some noteworthy works that suggest strategies to manage emergencies without explicitly dedicating OR time to them. \citet{lamiri2008} use a stochastic model to schedule elective patients in such a way to minimize OR overtime costs linked to emergencies. A more recent example is \citet{jha2023}, in which the authors propose the minimization of numerical indicator (TED, Total Expected Disturbance), calculated using historical data on emergencies and estimates of the probability of a given emergency happening in a certain moment. In this paper, we assume we {cannot} leverage on historical data to establish patterns or define probability of occurrence of emergencies, nor we can make assumptions on the amount of OR capacity an emergency can occupy. To represent the possible disruptions, we discretize the time at which emergencies can arise, and we assume that each emergency can be classified in the same moment it occurs into pre-defined duration categories (i.e. \textit{short}, \textit{medium} and \textit{long}). This allows to represent all the disruptions with a limited, and therefore enumerable, number of potential scenarios \citep{powell2007}. \\

For what concerns no-show, the topic is less studied in operations research literature than emergencies. \citet{dantas2018} provide a comprehensive review of no-show appointments in healthcare, without focusing on surgery only, to trace the main determinants that drive the phenomenon. Still not in the strict domain of surgical specialties, \citet{zacharias2014} propose a scheduling policy in which no-shows are counterbalanced by overbooking. \citet{berg2014} consider individually tailored no-show probabilities and lever on patient ordering, reducing the expected discomfort of patients by placing the patients that are more likely to result in no-shows towards the end of the patient sequence, if possible. As pointed out in \citet{ozen2016}, the availability of data on no-shows may be limited, because it is hard to determine if a patient will not show up for an appointment that has not happened yet, and they resort to simulation to test the impact of no-shows on their scheduling model. Some works (such as \citet{salazar2021}) proposed a novel classifier, leveraging on machine learning techniques, to predict no-shows for appointments, which, albeit very well performing, is hardly generalizable as is built with data coming from a specific non-surgical context. Moreover, as already mentioned above \citep{boudreau2011}, most same-day pediatric surgical cancellations (effectively equivalent to our definition of no-shows) are linked to acute illnesses, which are once again inherently hard to predict.\\
Without making strong assumptions on no-shows predictability, in our approach the number of possible no-shows corresponds to the number of patients scheduled in the first upcoming day(s) of the planning period, and therefore enumerable. We assume that each scenario is equally likely and find a way to proactively prepare a back-up plan, again allowing us to enumerate all potential disruptions of this second kind. \\
{Lastly, we would like to highlight that, to the best of our knowledge, no proposed approach incorporates the preparation of back-up scenarios in surgery scheduling.}

\section{Problem definition}
The problem is defined as the problem of scheduling elective surgeries from a waiting list over a planning horizon in a single ward of a Children's Hospital, without considering any upstream or downstream resources. This is done with the primary objective of minimizing an objective function made of penalty factors associated with how urgent the patients are and how much time they have spent in the waiting list. The more urgent a patient is, and the more time they have been waiting, the higher the penalty. Penalties are built in such a way that the more the assignment is further away in the future, the more they increase. The penalty associated with not scheduling a patient is always higher than the penalty associated with scheduling them on any day of the planning period. Overall, this penalty scheme promotes assignments closer in time, while granting equity to all patients in the list and without imposing hard constraints that may lead to infeasibilities. Surgery duration and OR capacity are considered when planning the schedule, as well as compatibility between days and OR and patients, due to the surgical teams assignment to rooms and days.\\

The schedule shall also be determined in such a way that it can be hastily re-deployed in case of disruptions. Two classes of disruptions are considered: \textit{emergencies}, defined as unexpected external patients that require immediate surgical treatment and that occupy OR time, and \textit{no-shows}, defined as scheduled patients that for any reason {cannot} undergo their surgery as planned. Back-up schedules need to be defined \textit{a priori} for each possible disruption, to allow for rescheduling on the fly. We consider disruptions happening on the first (current) day of the planning period. We assume that only one emergency or one no-show can happen on each given current day. This is consistent with the case of Buzzi hospital in Milan.\\
The creation of back-up schedules is {dependent} on the definition of a nominal schedule, that assigns patients to a day and an OR, and defines the order of patients assigned to the same day and OR, and, as a consequence, the starting and ending time of each surgery. Based on this, the occupation of the ORs is defined, and, with such plan, in every moment it is possible to determine if an OR is empty or which is its first availability if a patient is currently being operated on. \\
To model the problem, we assume that time is discretized, and therefore we are able to enumerate all the possible emergency realization scenarios, also by assuming that each emergency can have one of three pre-established lengths: short, medium, or long. The problem related to emergencies is double faceted: on the one side, to determine to which OR to assign each emergency that may happen, and, on the other, to determine how the schedule needs to be re-adjusted to accommodate for the disruption.\\
Similarly, all possible no-shows are readily listed, because each no-show scenario corresponds to an assignment done in the nominal schedule. To cope with no-shows, we introduce the concept of \textit{substitute} patient, a patient that is nominally assigned to the day after the one in which the emergency occurs, and that can be called in to fill the gap left by the no-show. The operative assumption done here is that the substitute patient is called in the hospital one day before their schedule surgery date. The no-show patient is then rescheduled after a pre-determined amount of days (or within a maximum number of days). Likewise, for each no-show scenario the rest of the schedule is restored to accommodate for the disruption.\\
Back-up schedules are designed so as to minimize the changes with respect to the nominal schedule, and to guarantee that patients assigned to the first days in the planning horizon are operated on within a limited number of days, even if they are affected by disruptions. Similarly, no-show patients must be rescheduled within a limited number of days. A limited OR overtime is accepted in all back-up schedules.

\section{MILP approach}
Let us now move on onto the formal description of the mathematical model of the problem. 
The primary aim is to assign patients (described by set $I$) to a nominal schedule over a planning horizon $D$, minimizing the penalty associated with (not) assigning them, weighted by their surgical priority and the time they spent in the list. Secondarily, it is necessary to design a feasible back-up schedule for each considered disruption, as close to the nominal one as possible.\\

Some information is given about each patient $i$ in the waiting list $I$. First, each patient's surgical deadline $l_i$ is known. It represents the clinically-recommended maximum amount of time it should elapse between being inserted into the list and the day of the surgery \citep{valente2009}. From this deadline, it is possible to determine each patient's urgency coefficient, $u_i$, equal to $360/l_i$ for each $i \in I$.
Secondly, it is known how many days patient $i$ has already spent in the waiting list, expressed by parameter $ m_i$. In order to determine the total waiting time for each patient, for each day $d$ of the planning horizon, it is possible to calculate a penalty factor associated with scheduling them on a given day $d$ (expressed with parameter $p_{id}$).  We also determine a penalty factor associated with not scheduling patient $i$ in the planning period, an information stored in parameter $q_i$:
\begin{align}
 & p_{id} = [d +\max{\{ m_i+d-l_i; 0}\}]u_i & \forall i \in I, d \in D \label{p} \\
 & q_{i} = [( m_i+|D|+1)+ \max\{( m_i+|D|+1-l_i);0\}]u_i & \forall i \in I \label{q}
\end{align}

As mentioned, the aim of the problem is to assign patients to a schedule in such a way that the penalty of (not) assigning them is as small as possible. We introduce binary decision variable $x_{idj}$, which equals 1 when patient $i$ is assigned to OR $j$ in day $d$, and formulate the following objective function, which was originally introduced in \citep{addis2016}:
\begin{align}
\min \sum_{i\in I}\sum_{j \in J} \left\{\sum_{d \in D} p_{id}x_{idj} +q_i\sum_{d \in D}\left(1-x_{idj}\right)\right\} \label{OF_}
\end{align}

Let us now move onto the description of constraints. For clarity, they will be divided into blocks. A handy summary of all sets and parameters used in the model is reported in Table \ref{tab:spvd}, together with the decision variables and their domains.\\

\textbf{Block \textit{A}} consists of patient assignment constraints:
\begin{align}
& \sum_{d \in D}\sum_{j \in J}x_{idj} \leq 1 & \forall i \in I \label{univocox}\\
& x_{idj} \leq a_{idj} & \forall i \in I, j \in J, d \in D \label{copertura}\\
& \sum_{i \in I}t_{idj}x_{idj}\leq T_{dj} & \forall j \in J, d \in D \label{capacitatot}
\end{align}

Constraints \eqref{univocox} make sure that each patient is scheduled at most once. Parameter $a_{idj}$ describes the compatibility between patients, days and ORs: it is equal to 1 if a surgical team can operate on patient $i$ on day $d$ and OR $j$. Constraints \eqref{copertura} ensure that only the proper patient-day-OR triples can be considered for assignment. Constraints \eqref{capacitatot} make sure that the sum of the surgical times $t_{idj}$ of all patients $i$ assigned to OR $j$ on day $d$ does not exceed total OR capacity $T_{dj}$. It is worth remarking that in this model time is not continuous, but rather it is represented by a set of time-slots, or \textit{moments}, $H$. All time-related parameters and decision variables used in this formulation and in the following are expressed as a function of $H$.\\

\textbf{Block \textit{B}} contains patient ordering constraints, which assign each scheduled patient to a position in the schedule. The possible positions are represented by a set $S$. $S$ is an ordered subset of $\mathbb{N}^+$, whose cardinality is big enough to provide a sufficient number of positions for any possible realization of the surgical plan. It is calculated as:
\begin{align}
 |S|=\left\lceil \max_{d \in D, j \in J}T_{dj}/\min_{i \in I, d \in D, j \in J}t_{idj} \right\rceil
\end{align}

Ordering is enforced only on a set $W\subset D$ of days, the first upcoming day(s) of the planning period, for which we want to guarantee back-up plans. We introduce three other decision variables. First, $y_{igjs}$ is a binary decision variable equal to 1 if patient $i$ assigned to OR $j$ on day $g$ is assigned to \textit{position} $s \in S$. Precedence variable $\nu_{irgj}$ is a binary variable equal to 1 if patient $i$ is assigned before patient $r$ in $g,j$. Based on the order, the starting time $\xi_{igj}$ of $i$'s surgery in $g,j$ can be computed. 
\begin{align}
&\sum_{s\in S} y_{igjs} = x_{igj} & \forall i \in I, j \in J, g \in W \label{corrispondenza}\\
&\sum_{i \in I}y_{igj(s+1)}\leq \sum_{i \in I}y_{igjs}& \forall j \in J, g \in W, s \in S\setminus \{|S|\} \label{primoposizionamento}\\
& \sum_{i\in I}y_{igjs}\leq 1 & \forall j \in J, g \in W, s \in S \label{univocoy}\\
& \xi_{igj} = \sum_{\substack{r\in I:}{r \neq i}}t_{rgj}\cdot \nu_{irgj} + x_{igj} & \forall i \in I, g \in W, j \in J \label{orarioinizio}\\
& y_{igjs} + \sum_{\substack{n \in S:}{n < s}}y_{rgjn} -1 \leq \nu_{irgj} & \forall i,r \in I: r \neq i, g \in W, j \in J, s \in S \label{nu1}\\
& \nu_{irgj} \leq x_{rgj} & \forall i,r \in I, g \in W, j \in J \label{nu2} \\
& \nu_{irgj} + \nu_{rigj} \leq 1 & \forall i,r \in I: r \neq i, g \in W, j \in J \label{nu3}\\
& \xi_{igj} \leq (T_{gj}+1-t_{igj})\cdot x_{igj} & \forall i \in I, g \in W, j \in J \label{completion_time}
\end{align}
Constraints \eqref{corrispondenza} link variables $y_{igjs}$ and $x_{igj}$, making sure that one and only one position is assigned to each scheduled patient. Constraints \eqref{primoposizionamento} guarantee that a position can be assigned only if the previous one is. Constraints \eqref{univocoy} ensure that each position is assigned at most once. Constraints \eqref{orarioinizio} calculate the starting time of each scheduled surgery. Constraints \eqref{nu1}-\eqref{nu3} calculate precedence variables $\nu_{irgj}$. Constraints \eqref{completion_time} force the starting time of $i$ in day $g$ and OR $j$ to be 0 if patient $i$ is not scheduled in $g,j$.\\

Blocks A and B determine the nominal schedule and its order. The following four blocks of constraints deal with the management of emergencies.\\

In order to manage an emergency, it is important to know what the current status of each OR {is} in the moment in which the emergency occurs. We define a non-negative variable $C_{gj}$ that encodes the slot $h$ in which the room will be released after the end of the last surgery. Complementarily, $\rho_{hgj}$ is a binary variable equal to 1 if OR $j$ on day $g$ is empty in time slot $h$. If the room is not empty, binary decision variables $\lambda_{igjh}$ are used to identify the patient $i$ that is currently being operated on in each OR $j$ at $h,g$. Knowing whether or not the room is empty, and, in the latter case, who is the patient currently being operated, allows to determine the quantity $ \Xi_{hgj}$, the first availability of room $j$ starting from moment $h$. To calculate the value of these variables, we introduce a series of constraints related to ORs, \textbf{Block \textit{C}}:
\begin{align}
& C_{gj} \geq \xi_{igj} + t_{igj}x_{igj} & \forall i \in I, j \in J, g \in W \label{isCi} \\
& h - C_{gj} \leq |H|\rho_{hgj} & \forall g \in W, h \in H, j \in J \label{isrho}\\
& h - C_{gj} \geq \rho_{hgj} - |H|(1-\rho_{hgj}) - 1 & \forall g \in W, h \in H, j \in J \label{isrho2}\\
& \sum_{i \in I} \lambda_{igjh} = 1-\rho_{hgj} & \forall g \in W, j \in J, h \in H \label{islambda}\\
& \xi_{igj} + t_{igj}x_{igj} \geq h \cdot \lambda_{igjh} & \forall i \in I, g \in W, j \in J, h \in H: h \neq 0 \label{rightbox}\\
& \xi_{igj} \leq h \cdot \lambda_{igjh} + |H|(1-\lambda_{igjh}) & \forall i \in I, g \in W, j \in J, h \in H: h \neq 0 \label{leftbox}\\
& \Xi_{gjh} \geq \xi_{igj} + t_{igj}x_{igj} - |H|(1-\lambda_{igjh}) & \forall i \in I, g \in W, j \in J, h \in H \label{whenXi}\\
& \Xi_{gjh} \leq \xi_{igj} + t_{igj}x_{igj} + |H|(1-\lambda_{igjh}) & \forall i \in I, g \in W, j \in J, h \in H \label{whenXi2}
\end{align}

Constraints \eqref{isCi}-\eqref{islambda} serve the purpose of calculating what is the time at which each OR becomes available after having completed all assigned surgeries, and of identifying which patients are being operated on in each OR $j$ and moment $h$, respectively. Constraints \eqref{rightbox} and \eqref{leftbox} identify, for each moment $h$, OR $j$ and day $g$, link the starting and ending time of patient $i$'s surgery with variable $\lambda_{igjh}$. Constraints \eqref{whenXi} and \eqref{whenXi2} calculate what is the first availability of OR $j$ in day $g$, starting from moment $h$.\\ 

An emergency that occurs in time slot $h$ of day $g$ must be managed and assigned to a OR: a binary variable $\eta_{hdj}$ is defined to determine which room the emergency that happens in moment $h$ is assigned to. The emergency assignment constraints of \textbf{Block \textit{D}} guarantee that each emergency is properly managed:
\begin{align}
& \sum_{j \in J}\eta_{hdj} = 1 & \forall h \in H, d \in W \label{emergenza}\\
& \Xi_{gjh} \leq \Xi_{gkh} + |H|(1-\eta_{hgj}) & \forall g \in W, j \in J, k \in J, h \in H \label{XiDef}
\end{align}

 Constraints \eqref{emergenza} force the assignment of each possible emergency to exactly one OR. Constraints \eqref{XiDef} force the emergency to be assigned to the first available OR.\\

\textbf{Block \textit{E}} contains constraints used to evaluate the impact of emergencies on scheduled patients. A new concept is introduced: emergency length class. As the length of an emergency {cannot} be known in advance, it is assumed to belong to a pre-determined class that approximates it. The set of the possible emergency lengths is defined as $L$. Together with $h$ (the moment the emergency occurs) and $g$ (the day), the triplet $h,g,l$ describes a specific scenario: an emergency of class $l \in L$ occurring in day $g$ at time-slot $h$. The impact is expressed using binary variable $\chi_{igj}$, which is equal to 1 if patient $i$ is impacted by the assignment of an emergency to OR $j$, that is, if the nominal assignment of that patient is scheduled to begin after $h$. The back-up schedules are described by binary assignment variables $\Bar{x}_{hglidj}$. For each scenario encoded by $h,g,l$, they store the information regarding its associated back-up schedule. Variables $\mu_{hgli}$ account for patients who {cannot} be rescheduled in the current planning period.
\begin{align}
& \sum_{j \in J}\xi_{igj} - h \leq |H|\cdot\chi_{hgi} & \forall i \in I, g \in W, h \in H \label{chi1}\\
& \sum_{j \in J}\xi_{igj} \geq h\cdot\chi_{hgi} + \varepsilon \cdot \chi_{hgi} & \forall i \in I, g \in W, h \in H \label{chi2}\\
& \sum_{\substack{d \in D:}{d < g+\Delta}}\sum_{j \in J} \Bar{x}_{hglidj} \geq \chi_{hgi} & \forall i \in I, h \in H, g \in W, l \in L \label{impattox}\\
& \sum_{\substack{k \in D:}{k \geq d; k \leq d+\Delta}} \sum_{j \in J} \Bar{x}_{hglikj} \geq \sum_{j \in J}x_{idj} - \mu_{hgli}& \forall h \in H, g \in W, i \in I, d \in D, l \in L \label{impattoxbar}
\end{align}

Constraints \eqref{chi1} and \eqref{chi2} force the suitable value of variables $\chi_{igj}$. In \eqref{chi2}, $\varepsilon$ is a sufficiently small constant. Constraints \eqref{impattox} force each impacted patient to be considered in the rescheduling, while Constraints \eqref{impattoxbar} determine if and how many patients are excluded ($\mu_{hgli}$) from the rescheduling.\\

Then, we introduce \textbf{Block \textit{F}}, which contains the constraints related to the rescheduling caused by the emergency:
\begin{align}
& \sum_{j \in J}x_{idj} + \sum_{j \in J}\sum_{k \in D: k < d}\Bar{x}_{hglikj} \leq 1 & \forall i \in I, d \in D, h \in H, g \in W, l \in L \label{nonanticipare}\\
& \sum_{j \in J}\sum_{k \in D}\Bar{x}_{hglikj} \leq \sum_{j \in J}\sum_{d \in D}x_{idj} & \forall h \in H, g \in W, l \in L, i \in I \label{nonmetterlo} \\
& \Bar{x}_{hgligj} \geq x_{igj}-\chi_{hgi} & \forall h \in H, g \in W, i \in I, j \in J, l \in L \label{mantienixbar}\\
& \sum_{d \in D}\sum_{j \in J} \Bar{x}_{hglidj} \leq 1 & \forall h \in H, g \in W, i \in I, l \in L \label{univocoxbar}\\
& \Bar{x}_{hglidj} \leq a_{idj} & \forall h \in H, g \in W, i \in I, j \in J, d \in D, l \in L \label{coperturaxbar}\\
& \sum_{i \in I} t_{idj}\Bar{x}_{hglidj} \leq T_{dj} & \forall j \in J, d \in D, h \in H, g \in W, l \in L: d > g \label{capacitatotxbar}\\
& \sum_{i \in I} t_{igj}\Bar{x}_{hgligj} \leq T_{gj} - \min{\{\gamma_l; (T_{gj} - h)\}}\cdot\eta_{hgj} + \Omega & \forall j \in J, d \in D, h \in H, g \in W, l \in L \label{capacitatot1xbar}
\end{align}

Constraints \eqref{nonanticipare} prevent assignments to be anticipated with respect to their original date. Constraints \eqref{nonmetterlo} prevent patients that were not previously included in the schedule to be included in the new schedule. Constraints \eqref{mantienixbar} maintain the original assignment of patients that have not been affected by the emergency considered in the scenario. Constraints \eqref{univocoxbar}-\eqref{capacitatotxbar} are the corresponding of Constraints \eqref{univocox}-\eqref{capacitatot} of the nominal assignment problem, while Constraints \eqref{capacitatot1xbar} allow for some overtime (defined by parameter $\Omega$) to be accepted in the day of the emergency, as the emergency itself occupies some OR time. This was the last block of Constraints concerning emergency management, the following blocks dealing instead with the management of no-shows.\\

The number of possible no-show scenario is enumerable by definition. Only patients that have been assigned in the nominal schedule in $g\in W$ can become no-shows, so, differently from what we did for the modelling of emergency constraints, it is not necessary to make any additional assumptions. The back-up schedules are described by binary assignment variables $\hat{x}_{bgidj}$. For each scenario encoded by no-show patient $b$ in day $g$, they store the information regarding its associated back-up schedule. \textbf{Block \textit{G}} defines a substitute patient for all rooms, which is identified by binary variable $\theta_{igj}$:
\begin{align}
& \sum_{i \in I}\theta_{ijg} = 1 & \forall j \in J, g \in W \label{nomultitheta}\\
& \sum_{j \in J} \theta_{ijg} \leq \sum_{k \in J} x_{i(g+1)k} & \forall i \in I, g \in W \label{impattotheta}\\
& \theta_{ijg} \leq a_{ijg} & \forall g \in W, j \in J, i \in I \label{coperturatheta} \\
& \hat{x}_{bgigj} \geq \theta_{igj} + x_{bgj} - 1 & \forall b \in I, g \in W, i \in I, j \in J \label{agganciotheta} \\
& \sum_{j \in J} \hat{x}_{bgb(g+\hat{\Delta})j} \geq \sum_{j \in J} x_{bgj} & \forall b \in I, g \in W, i \in I \label{force_delay}
\end{align}

Constraints \eqref{nomultitheta} force the identification of a substitute patient for each room to be unique. Constraints \eqref{impattotheta} force the substitute patient to be a patient originally scheduled in the following day with respect to day $g\in W$. Constraints \eqref{coperturatheta} ensure that the substitute patient is compatible with the OR they may be called to cover. Constraints \eqref{agganciotheta} make sure that the substitute patient is scheduled instead of the patient they are called to cover for, in the scenario in which the latter is a no-show. Constraints \eqref{force_delay} ensure that the no-show patient is rescheduled after $\hat{\Delta}$ days. \eqref{force_delay} can be substituted with the following \eqref{force_delay_fixed}, in case we do not want to fix the delay, but just ensure that the no-show patient is rescheduled within $\Delta$ days:
\begin{align}
& \sum_{j \in J}\sum_{d \in D: d < g + \Delta} \hat{x}_{bgbdj} \geq \sum_{j \in J} x_{bgj} & \forall b \in I, g \in W, i \in I \label{force_delay_fixed} 
\end{align}

\textbf{Block \textit{H}} deals with back-up plans in case of reassignment whenever a no-show happens:
\begin{align}
& \sum_{d \in D}\sum_{j \in J} \hat{x}_{bgidj} \leq 1 & b \in I, g \in W, i \in I \label{univocoxhat}\\
& \hat{x}_{bgidj} \leq a_{idj} & \forall h \in H, g \in W, i \in I, j \in J, d \in D \label{coperturaxhat} \\
& \sum_{i \in I} t_{idj}\hat{x}_{bgidj} \leq T_{dj} & \forall j \in J, d \in D, h \in H, g \in W: d > g \label{capacitaxhat}\\
& \sum_{i \in I} t_{igj}\hat{x}_{bgigj} \leq T_{gj} + \Omega & \forall j \in J, h \in H, g \in W \label{capacitaxhat2}\\
& \sum_{j \in J}x_{idj} + \sum_{j \in J}\sum_{k \in D: k < d}\hat{x}_{bgikj} \leq 1 & \forall i \in I, d \in D, b \in I, g \in W \label{nonanticiparexhat}\\
& \sum_{j \in J}\sum_{k \in D}\hat{x}_{bgikj} \leq \sum_{j \in J}\sum_{d \in D}x_{idj} & \forall g \in W, b \in I, i \in I \label{nonmetterlotheta} \\
& \hat{x}_{bgigj} \geq x_{igj} & \forall b \in I, g \in W, i \in I, j \in J: i \neq b \label{mantienixhat}
\end{align}

Constraints \eqref{univocoxhat}-\eqref{capacitaxhat} are the corresponding of Constraints \eqref{univocox}-\eqref{capacitatot}. Constraints \eqref{capacitaxhat2} allow for some overtime to accommodate for the substitute patient, in a similar fashion to what was done in \eqref{capacitatot1xbar}. Constraints \eqref{nonanticiparexhat} prevent patients' admission dates in the reschedule to be anticipated with respect to the original plan. Constraints \eqref{nonmetterlotheta} prevent patients that were not previously included in the nominal schedule to be included in the new schedule. Lastly, Constraints \eqref{mantienixhat} ensure that all patients included in the nominal scheduled are included in the reschedule.\\

\begin{table}[p]
\centering
\small
\begin{tabular}{p{0.05\textwidth}p{0.92\textwidth}}
\toprule
\textbf{Sets} & \\
\midrule
$I$ & patients in the waiting list \\
$J$ & ORs \\
$D$ & days (indexed by $d$) of the planning horizon \\
$W$ & days (indexed by $g$) in which we want to enforce robustness \\
$S$ & position in the daily plan occupied by a patient\\
$H$ & time-slot discretization of the duration of the day. Starts at 0 \\
$L$ & set of possible types of emergency lengths \\
\midrule
\multicolumn{2}{l}{\textbf{Parameters}} \\
\midrule
$a_{idj}$ & compatibility matrix between patient $i$, day $d$ and OR $j$\\
$ m_i$ & days spent in the waiting list by patient $i$\\
$u_i$ & priority class of patient $i$\\
$l_i$ & surgery deadline of patient $i$\\
$q_i$ & penalty for not assigning patient $i$\\
$p_{id}$ & penalty for assigning patient $i$ to day $d$ \\
$T_{dj}$ & availability of OR $j$ in day $d$, expressed in time-slots\\
${t}_{idj}$ & duration of $i$'s surgery in day $d$ and OR $j$, expressed in time-slots\\
$\gamma_l$ & duration of emergency in scenario $l$\\
$\Omega$ & acceptable overtime in case of a substitution of a no-show or an emergency, expressed in time-slots\\
$\hat{\Delta}$ & fixed delay used to reschedule no-show patients, expressed in days \\
$\Delta$ & maximum rescheduling delay\\
\midrule
\multicolumn{2}{l}{\textbf{Decision variables}} \\
\midrule
$x_{idj}$ & binary, assignment of patient $i$ to day $d$ and OR $j$\\
$y_{igjs}$ & binary, assignment of patient $i$ in day $g$ to OR $j$ to slot $s$ (ordering) \\
$\nu_{irgj}$ & binary, equal to 1 if patient $i$ is assigned before patient $r$ in day $g$ and OR $j$\\
${\xi}_{igj}$ & continuous non-negative, expected starting time (in time-slots) of the surgery for $i$ in $g,j$ \\
$C_{gj}$ & continuous non-negative, end time of last patient in room $j$ in day $g$ \\
$\rho_{hgj}$ & binary, equal to 1 if room $j$ is empty when an emergency in time slot $h$ of day $g$ occurs \\
$\lambda_{igjh}$ & binary, equal to 1 if patient $i$ is currently being operated on in OR $j$ when an emergency occurs in time-slot $h$ of day $g$\\
$ \Xi_{hgj}$ & continuous non-negative, first availability in room $j$ to assign emergency occurring in time-slot $h$ of day $g$ \\
$\eta_{hgj}$ & binary, equal to 1 if an emergency occurring in time slot $h$ of day $g$ is assigned to $j$\\
$\chi_{hgi}$ & binary, equal to 1 if the assignment of $i$ is affected by emergency occurring in time slot $h$ of day $g$\\
$\Bar{x}_{hglidj}$ & binary, equal to 1 if $i$ needs to be rescheduled in day $d$ and OR $j$ because of emergency occurring in time slot $h$ of day $g$\\
$\mu_{hgli}$ & binary, equals 1 if patient $i$ {cannot} be rescheduled if an emergency of type $l$ occurs in time slot $h$ of day $g$\\
$\hat{x}_{bgidj}$ & binary, equal to 1 if $i$ needs to be rescheduled in day $d$ and OR $j$ because of a no-show patient $b$ in day $g$\\
$\theta_{igj}$ & binary, equal to 1 if patient $i$ is the substitute patient for day $g$ and OR $j$ in case of no-shows\\
\bottomrule
\end{tabular}
\caption{Summary of all sets, parameters and decision variables used in the MILP model.}
\label{tab:spvd}
\end{table}

\subsection{Solution approach}
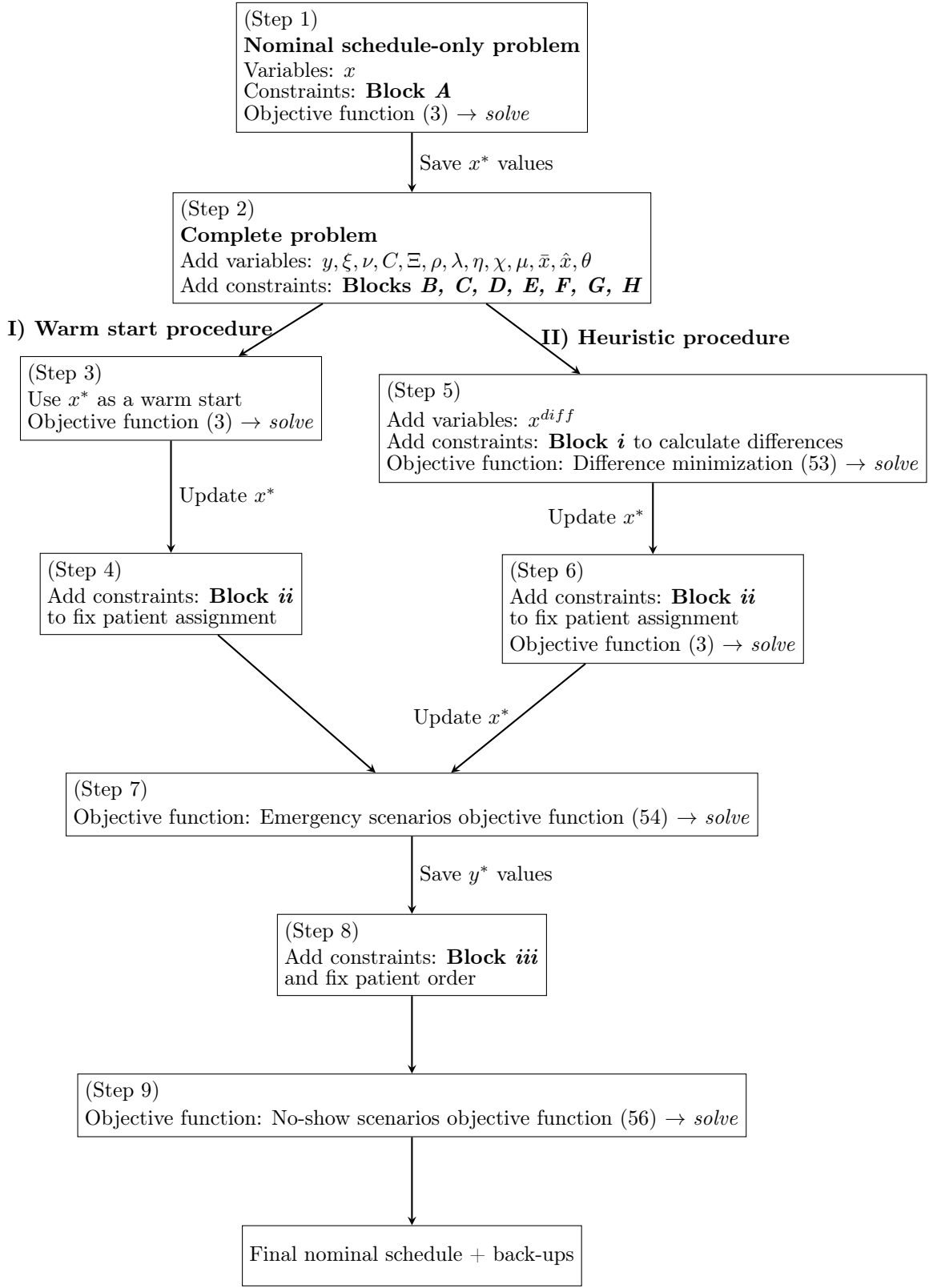
\begin{figure}[htbp]
\small
{\centering
\begin{tikzpicture}[node distance=2.5cm]
\node (nominal) [process] {\shortstack[l]{(Step 1) \\\textbf{Nominal schedule-only problem} \\ Variables: $x$ \\ Constraints: \textbf{Block \textit{A}} \\ Objective function \eqref{OF_} $\rightarrow$ \textit{{solve}}}};

\node (full) [process, below of=nominal,yshift=-0.5cm] {\shortstack[l]{(Step 2) \\\textbf{Complete problem} \\ Add variables: $y, \xi, \nu, C, 
 \Xi, \rho, \lambda, \eta, \chi, \mu, \Bar{x}, \hat{x}, \theta$ \\ Add constraints: \textbf{Blocks \textit{B, C, D, E, F, G, H}}}};
 
\node (fixing1) [process, below of=full, xshift = 4cm, yshift = -.5cm] {\shortstack[l]{(Step 5) \\Add variables: $x^{diff}$ \\ Add constraints: \textbf{Block \textit{i}} to calculate differences \\ Objective function: Difference minimization \eqref{OF_heu} $\rightarrow$ \textit{{solve}}}};

\node (fixing15) [process, below of=fixing1,yshift = -.5cm] {\shortstack[l]{(Step 6) \\Add constraints: \textbf{Block \textit{ii}}\\to fix patient assignment\\ Objective function \eqref{OF_} $\rightarrow$ \textit{solve}}};

\node (ws) [process, below of=full, xshift = -4cm] {\shortstack[l]{(Step 3) \\Use $x^*$ as a warm start\\ Objective function \eqref{OF_} $\rightarrow$ \textit{solve}}};

\node (fixing2) [process, below of=ws, yshift=-0.75cm] {\shortstack[l]{(Step 4) \\Add constraints: \textbf{Block \textit{ii}} \\ to fix patient assignment}};

\node (xbar) [process, below of=fixing2, xshift =4cm, yshift = -1cm] {\shortstack[l]{(Step 7) \\Objective function: Emergency scenarios objective function \eqref{OF_xbar} $\rightarrow$ \textit{{solve}}}};

\node (fixing3) [process, below of=xbar] {\shortstack[l]{(Step 8) \\ Add constraints: \textbf{Block \textit{iii}}\\ and fix patient order}};

\node (xhat) [process, below of =fixing3] {\shortstack[l]{(Step 9) \\Objective function: No-show scenarios objective function \eqref{OF_xhat} $\rightarrow$ \textit{{solve}}}};

\node (final) [process,below of=xhat] {Final nominal schedule + back-ups};
\draw [arrow] (nominal) -- node[anchor=west] {Save $x^*$ values} (full);
\draw [arrow] (full) -- node[anchor=west] {\textbf{II) Heuristic procedure}} (fixing1);
\draw [arrow] (full) -- node[anchor=east]{\textbf{I) Warm start procedure}}(ws);
\draw [arrow] (fixing1) -- node[anchor=east] {Update $x^*$} (fixing15);
\draw [arrow] (fixing15) -- node[anchor=east] {Update $x^*$} (xbar);
\draw [arrow] (ws) -- node[anchor=west] {Update $x^*$} (fixing2);
\draw [arrow] (fixing2) -- node[anchor=east] {}(xbar);
\draw [arrow] (xbar) -- node[anchor=west]{Save $y^*$ values}(fixing3);
\draw [arrow] (fixing3) -- node[anchor=east]{}(xhat);
\draw [arrow] (xhat) -- node[anchor=east]{}(final);
\end{tikzpicture}\\}
\caption{Schematic representation of the heuristic algorithm.}
\label{img:heu2}
\end{figure}

The model {cannot} be solved as a whole but in some smaller instances. However, even rather small optimality gaps lead to solutions that are unacceptable for decision-makers, as the proposed plan overuses spare capacity to guarantee a feasible rescheduling for the days in which a back-up plan is required. Consequently, we developed a heuristic approach that first computes a nominal schedule and then uses such schedule to compute a feasible solution for the overall problem applying:
\begin{enumerate}[I)]
 \item a warm start procedure, or
 \item a MILP-based sequential heuristic.
\end{enumerate}
The back-up schedules are then optimized based on patients' penalties. Table \ref{tab:spvd_heu} provides a summary of the parameters and decision variables used by the heuristic, while Figure \ref{img:heu2} represents graphically the proposed approach.\\
In details, we first solve the \textit{nominal} problem, defined by {Block \textit{A}} and objective function \eqref{OF_} (Step 1). The solution of this problem is stored in parameter $x_{idj}^*$. Then, the warm start and the heuristic are applied to obtain a feasible complete solution.\\

The warm start computes a feasible solution of the \textit{complete} problem (Step 3), whose assignments are then (Step 4) set through the following additional constraints (\textbf{Block \textit{ii}}):
\begin{align}
 & \sum_{j \in J}x_{igj} \geq \sum_{j \in J} x^*_{igj} & \forall i \in I, g \in W \label{fix_xx} 
\end{align}
Constraints \eqref{fix_xx} fix the admission dates of the patients included in the robust days of the nominal solution of the warm start procedure.\\

The heuristic procedure introduces variables $x^{diff}_{id}$ to represent the differences between the \textit{nominal} schedule (\textbf{Block A}-only problem) and any feasible solution of the \textit{complete} one (\textbf{Blocks A, B, C, D, E, F, G, H}). The differences are computed through the following \textbf{Block \textbf{\textit{i}}} of constraints:
\begin{align}
 & \sum_{j \in J}x_{idj} - \sum_{j \in J}x^*_{idj} \leq x^{diff}_{id} & \forall i \in I, d \in D \label{differenzax1} \\
 & \sum_{j \in J}x^*_{idj} - \sum_{j \in J}x_{idj} \leq x^{diff}_{id} & \forall i \in I, d \in D \label{differenzax2}
\end{align}

In Step 5, a new problem is considered. It includes all the constraints of Blocks A-H and the additional constraints of Block \textit{ii}. Its goal is to minimize at the same time: \textit{1}) the difference between the two solutions and \textit{2}) the number of patients that {cannot} be rescheduled, with the following objective function:
\begin{align}
 \min \sum_{i \in I}\sum_{d \in D}x^{diff}_{id}\label{OF_heu}
\end{align}

The solution of this new \textit{complete} schedule is stored, updating parameters $x_{idj}^*$. 
{Block \textit{ii}} of constraints is added and the problem is solved once more with objective function \eqref{OF_}. Here, the greater or equal to sign of \eqref{fix_xx} allows for the model to try and find solutions in which, at the cost of introducing differences with respect to the initial solution, OR capacity is used more. The value of $x_{igj}^*$ is updated once again (Step 6). \\

After having fixed a common starting point, the \textit{final} nominal schedule, and regardless of whether it was obtained with the warm start or the heuristic procedure, we can now optimize for each type of disruption. First, we optimize for the emergency scenarios and then for the no-shows. This is a deliberate choice, because to better manage the emergencies it is necessary to establish the ordering $y_{idjs}$ of patients, a decision that does not impact on no-show management.\\
To optimize for emergencies, we solve the complete model with the following objective function (Step 7), defined as a sum of penalties related to (not) assigning patients in emergency back-up schedules, in a similar fashion to what was done in \eqref{OF_}:
\begin{align}
\min \frac{1}{|H|\cdot|W|\cdot|L|}\sum_{i\in I}\sum_{g \in W}\sum_{h \in H}\sum_{j \in J}\left\{\sum_{d \in D: d > g} p_{id}\bar{x}_{hgidj} +q_i\sum_{d \in D}\left(1-\bar{x}_{hgidj}\right)\right\} \label{OF_xbar}
\end{align}

As mentioned, the outcome of Step 7 also orders patients within their assigned days, a decision stored in parameter $y^*_{igjs}$. \textbf{Block \textit{iii}} fixes the order and the OR assignment of patients (Step 8) with the following additional constraints:
\begin{align}
 & y_{igjs} \geq y^*_{igjs} & \forall i \in I, g \in W, j \in J, s \in S \label{fix_yy} 
\end{align}

To optimize for no-shows, we solve the complete model with the following objective function (Step 9):
 \begin{align}
 \min \frac{1}{|I|/|D|}\sum_{i\in I}\sum_{g \in W}\sum_{b \in I}\sum_{j \in J}\left\{\sum_{d \in D} p_{id}\hat{x}_{bgidj} +q_i\sum_{d \in D}\left(1-\hat{x}_{bgidj}\right)\right\} \label{OF_xhat}
 \end{align}

\begin{table}[h!]
\centering
\small
\begin{tabular}{p{0.07\textwidth}p{0.88\textwidth}}
\toprule
\multicolumn{2}{l}{\textbf{Parameters}} \\
\midrule
$x_{idj}^*$ & Binary parameter storing either the nominal schedule or {the} \textit{final} nominal schedule \\
$y_{igjs}^*$ & Binary parameter storing the order of the \textit{final} nominal schedule \\
\midrule
\multicolumn{2}{l}{\textbf{Decision variables}} \\
\midrule
$x^{diff}_{id}$ & {B}inary, equal to 1 if the assignment of patient $i$ in day $d$ differs from the nominal-only solution \\
\bottomrule
\end{tabular}
\caption{Parameters and decision variables used in heuristic algorithm.}
\label{tab:spvd_heu}
\end{table}

\section{Results}
\begin{figure}[h]
 \centering
\begin{subfigure}{.4\textwidth}
\centering
\includegraphics[width=.8\textwidth]{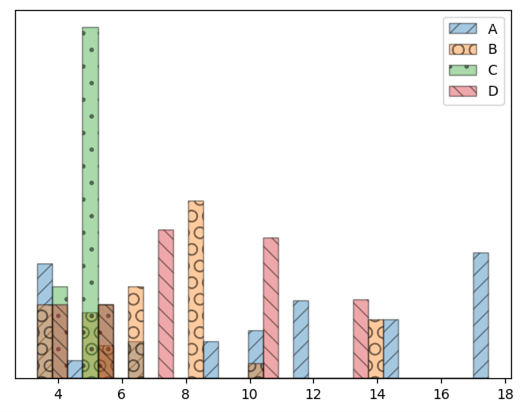} 
\caption{Distributions of surgical times for the four patient groups, expressed as multiples of 15 minutes.}
 \label{fig:hist}
\end{subfigure}
\hspace{3em}
\begin{subfigure}{.4\textwidth}
\centering
\includegraphics[width=.6\textwidth]{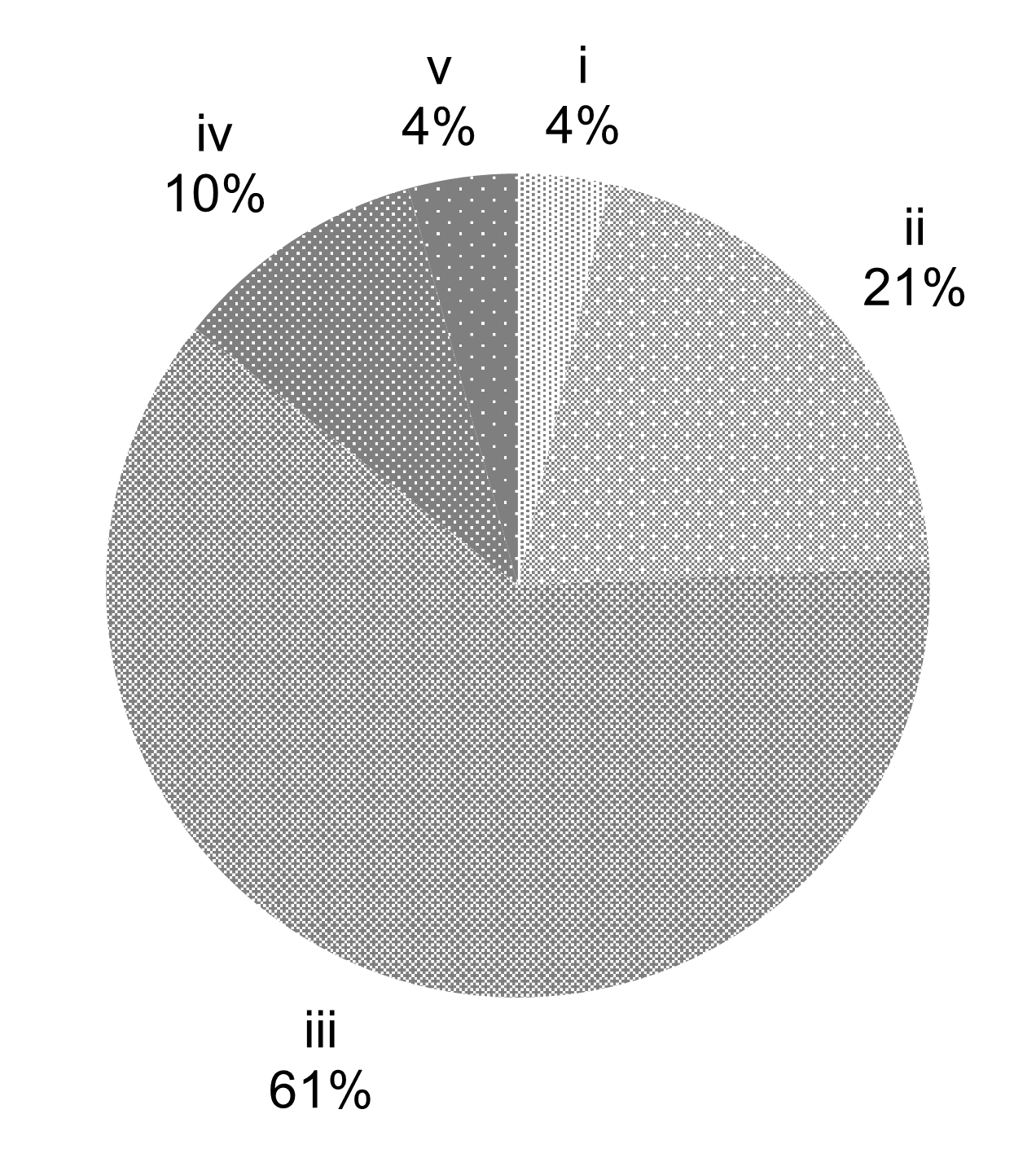} 
\caption{Distributions of the different urgency classes in the tested instances, \textit{i} being the most urgent.}
 \label{fig:pie}
\end{subfigure}
\label{fig:test}
\caption{Characterization of the patients in the test instances.}
\end{figure}

We first tested our approaches on a set of instances, based on \citet{addis2016}. In this instances, patients belong to one of five urgency classes as proposed in \citet{valente2009}, which are true to the Italian guidelines on surgical priorities.\\
In the instances, $|D|$ is equal to $[14,28]$, $|J|=[2,3]$, $|W| =1$ for all combinations. $a_{idj}$ was always set to 1 but in those days indexed by a $d$ multiple of 6 or 7, to model the weekends. The time discretization (set $H$) is assumed to be expressed in terms of 15 minutes blocks: as we consider 6 hours of operating rooms activity, we set $T_{dj}=24$; as one hour of overtime can be accepted for rescheduling, $\Omega$ was set to $4$. $\Delta$ and $\hat{\Delta}$ were set to $7$ and $2$, respectively, meaning that no-show patients must be rescheduled in two days, and patients cancelled because of an emergency must be rescheduled within a week. Lastly, the possible duration of emergencies were defined as $\gamma_l = [4,8,16]$. We considered instances with 40, 80 and 120 patients. For each cardinality of $I$, we consider four distinct patient populations, each characterized by different patient mixes in terms of surgical length, yielding to four types of instances, A, B, C and D, as shown in Figure \ref{fig:hist}. \ref{fig:pie} reports the distribution of the urgency classes in all the tested instances, which is the same for all groups of patients.

Approaches and models were implemented in \texttt{gurobipy}~\cite{gurobipy-web} and solved using Gurobi 11.0.1, in a \texttt{jupyter} environment. Experiments were run on a Microsoft Windows 11 machine with an Intel\textsuperscript{\textregistered} Core\textsuperscript{TM} i5-1135G7 CPU and 16 GB of installed RAM. For all steps of the approach, a time limit of 900 seconds was imposed, without enforcing any memory limit. We compared the results of the heuristic procedure with those obtained with the warm start procedure. Results are reported in Tables \ref{tab:gaps40}, \ref{tab:gaps80} and \ref{tab:gaps120}.\\

The Tables are structured as follows: the ID column reports the ID of the tested patient list. This, along with column $|D|$ and $|J|$, defines the tested instance. $OF$ is the objective function of the nominal schedule problem, while $OF_{WS}$ and $LB_{WS}$ represent the objective function value of the complete problem (as calculated by the warm start procedure) and its lower bound, respectively. Similarly, $OF_{Heu}$ is the objective function value of the complete problem calculated by the heuristic procedure. Columns $OF_{WS}-OF$ and $OF_{Heu}-OF$ calculate the difference between the two objective functions of the complete problem and that of the nominal schedule, in some sense calculating the \textit{"price of robustness"} \cite{bertsimas2004price} of introducing a number of constraints to guarantee that a feasible back-up schedule exists for each disruption. 
Lastly, columns $\frac{OF_{WS}-LB_{WS}}{LB_{WS}}$ and $\frac{OF_{Heu}-LB_{WS}}{LB_{WS}}$ report the percentage gap between the objective function values computed by the warm start and the heuristic procedure and the lower bound found by the former.\\

In all tested instances, the nominal schedule problem (Step 1) is solved to optimality in an average of $21.5$ seconds but in three instances. In these, the solver reaches the time limit, with an average gap of $0.02\%$. Including those three instances, all with cardinality $|I|=120$, the average computational time to solve the nominal problem is $76.4$ seconds. Still considering all instances at once, the warm start procedure (Step 3) always reaches the time limit with an average gap of 7.9\% (maximum 135.8\%, minimum 0.07\%), while the heuristic procedure takes on average of $732.7$ seconds to solve. The first part of the heuristic (Step 5) is solved to optimality in 27 instances out of 48.\\

Results for $|I|=40$ are reported in Table \ref{tab:gaps40}.
The warm start procedure reaches a lower objective function with respect to the heuristic procedure in 11 instances out of 16, with an average gap from the lower bound of the complete problem of $3.72\%$ versus an average gap of $3.99\%$ reached by the heuristic. The difference between the two approaches is small, and judging on these first results only would not establish a clear winner. \\

\begin{table}[h]
\centering
\scriptsize
\begin{tabular}{cccrrrrrrrrr}
\hline
\multicolumn{1}{l}{\textbf{ID}} & \multicolumn{1}{c}{\textbf{$|D|$}} & \multicolumn{1}{c}{\textbf{$|J|$}} & \multicolumn{1}{l}{\textbf{$OF$}} & \multicolumn{1}{c}{\textbf{$LB_{WS}$}} & \multicolumn{1}{c}{\textbf{$OF_{WS}$}} & \multicolumn{1}{c}{\textbf{$OF_{Heu}$}} & \multicolumn{1}{c}{\textbf{$OF_{WS}-OF$}} & \multicolumn{1}{c}{\textbf{$OF_{Heu}-OF$}} & \textbf{$\frac{OF_{WS}-LB_{WS}}{LB_{WS}}$} & \textbf{$\frac{OF_{Heu}-LB_{WS}}{LB_{WS}}$} \\
\hline
\textbf{40 A} & 14 & 2 & 4891 & 4886 & 4913 & 4939 & 22 & 48 & 0.55\% & 1.08\% \\
\textbf{40 A} & 14 & 3 & 4584 & 4579 & 4611 & 4668 & 27 & 84 & 0.70\% & 1.94\% \\
\textbf{40 A} & 28 & 2 & 5067 & 5059 & 5096 & 5151 & 29 & 84 & 0.73\% & 1.82\% \\
\textbf{40 A} & 28 & 3 & 4659 & 4651 & 4783 & 4707 & 124 & 48 & 2.84\% & 1.20\% \\
\textbf{40 B} & 14 & 2 & 3200 & 3200 & 3223 & 3387 & 23 & 187 & 0.72\% & 5.84\% \\
\textbf{40 B} & 14 & 3 & 3004 & 3003 & 3332 & 3202 & 328 & 198 & 10.96\% & 6.63\% \\
\textbf{40 B} & 28 & 2 & 3289 & 3286 & 3323 & 3938 & 34 & 649 & 1.13\% & 19.84\% \\
\textbf{40 B} & 28 & 3 & 3033 & 3028 & 3689 & 3348 & 656 & 315 & 21.83\% & 10.57\% \\
\textbf{40 C} & 14 & 2 & 3138 & 3138 & 3156 & 3166 & 18 & 28 & 0.57\% & 0.89\% \\
\textbf{40 C} & 14 & 3 & 3005 & 3005 & 3015 & 3045 & 10 & 40 & 0.33\% & 1.33\% \\
\textbf{40 C} & 28 & 2 & 3214 & 3214 & 3232 & 3286 & 18 & 72 & 0.56\% & 2.24\% \\
\textbf{40 C} & 28 & 3 & 3024 & 3024 & 3300 & 3060 & 276 & 36 & 9.13\% & 1.19\% \\
\textbf{40 D} & 14 & 2 & 5768 & 5768 & 5773 & 5852 & 5 & 84 & 0.09\% & 1.46\% \\
\textbf{40 D} & 14 & 3 & 5486 & 5481 & 5506 & 5638 & 20 & 152 & 0.46\% & 2.86\% \\
\textbf{40 D} & 28 & 2 & 5863 & 5863 & 5899 & 5971 & 36 & 108 & 0.61\% & 1.84\% \\
\textbf{40 D} & 28 & 3 & 5523 & 5522 & 5978 & 5691 & 455 & 168 & 8.26\% & 3.06\% \\
\multicolumn{3}{c}{\textbf{Average}} & \textbf{4172} & \textbf{4169} & \textbf{4302} & \textbf{4316} & \textbf{130} & \textbf{144} & \textbf{3.72\%} & \textbf{3.99\%} \\
\hline
\end{tabular}
\caption{Results for $|I|=40$.}
\label{tab:gaps40}
\end{table}

For $|I|=80$, the warm start procedures still performs better than the heuristic one in 11 cases out of 16, but the average gap with respect to the lower bound of the complete problem is higher for the warm start ($4.97\%$) with respect to the heuristic ($4.05\%$). This is due to the fact that the highest gaps of the warm start are larger than those of the heuristic.
For instance, the maximum gap for the warm start is $17.97\%$, versus a maximum of $8.95\%$ (in different instances) of the heuristic. This superiority is also evident when looking at the average values of the \textit{price of robustness} columns ($OF_{WS}-OF$ and $OF_{Heu}-OF$): on average, the heuristics produces solutions with less increment in the objective function.\\

\begin{table}[h]
\centering
\scriptsize
\begin{tabular}{cccrrrrrrrrr}
\hline
\multicolumn{1}{l}{\textbf{ID}} & \multicolumn{1}{c}{\textbf{$|D|$}} & \multicolumn{1}{c}{\textbf{$|J|$}} & \multicolumn{1}{l}{\textbf{$OF$}} & \multicolumn{1}{c}{\textbf{$LB_{WS}$}} & \multicolumn{1}{c}{\textbf{$OF_{WS}$}} & \multicolumn{1}{c}{\textbf{$OF_{Heu}$}} & \multicolumn{1}{c}{\textbf{$OF_{WS}-OF$}} & \multicolumn{1}{c}{\textbf{$OF_{Heu}-OF$}} & \textbf{$\frac{OF_{WS}-LB_{WS}}{LB_{WS}}$} & \textbf{$\frac{OF_{Heu}-LB_{WS}}{LB_{WS}}$} \\
\hline
\textbf{80 A} & 14 & 2 & 9677 & 9677 & 9716 & 9886 & 39 & 209 & 0.40\% & 2.16\% \\
\textbf{80 A} & 14 & 3 & 8888 & 8881 & 8965 & 9383 & 77 & 495 & 0.95\% & 5.65\% \\
\textbf{80 A} & 28 & 2 & 10631 & 10629 & 10682 & 10807 & 51 & 176 & 0.50\% & 1.67\% \\
\textbf{80 A} & 28 & 3 & 9428 & 9420 & 10310 & 9550 & 882 & 122 & 9.45\% & 1.38\% \\
\textbf{80 B} & 14 & 2 & 8947 & 8939 & 9086 & 9156 & 139 & 209 & 1.64\% & 2.43\% \\
\textbf{80 B} & 14 & 3 & 8315 & 8295 & 9219 & 8507 & 904 & 192 & 11.14\% & 2.56\% \\
\textbf{80 B} & 28 & 2 & 9494 & 9485 & 9602 & 9784 & 108 & 290 & 1.23\% & 3.15\% \\
\textbf{80 B} & 28 & 3 & 8637 & 8615 & 9712 & 9052 & 1075 & 415 & 12.73\% & 5.07\% \\
\textbf{80 C} & 14 & 2 & 8429 & 8428 & 8498 & 8621 & 69 & 192 & 0.83\% & 2.29\% \\
\textbf{80 C} & 14 & 3 & 7845 & 7825 & 9231 & 7929 & 1386 & 84 & 17.97\% & 1.33\% \\
\textbf{80 C} & 28 & 2 & 8832 & 8823 & 8920 & 9048 & 88 & 216 & 1.10\% & 2.55\% \\
\textbf{80 C} & 28 & 3 & 8033 & 7957 & 8189 & 8669 & 156 & 636 & 2.92\% & 8.95\% \\
\textbf{80 D} & 14 & 2 & 9405 & 9394 & 9408 & 9828 & 3 & 423 & 0.15\% & 4.62\% \\
\textbf{80 D} & 14 & 3 & 8730 & 8713 & 8959 & 9121 & 229 & 391 & 2.82\% & 4.68\% \\
\textbf{80 D} & 28 & 2 & 10184 & 10176 & 10523 & 10514 & 339 & 330 & 3.41\% & 3.32\% \\
\textbf{80 D} & 28 & 3 & 9153 & 9112 & 10226 & 10296 & 1073 & 1143 & 12.23\% & 12.99\% \\
\multicolumn{3}{c}{\textbf{Average}} & \textbf{9039} & \textbf{9023} & \textbf{9453} & \textbf{9384} & \textbf{414} & \textbf{345} & \textbf{4.97\%} & \textbf{4.05\%} \\
\hline
\end{tabular}
\caption{Results for $|I|=80$.}
\label{tab:gaps80}
\end{table}

For $|I|=120$, the warm start procedures performs better than the heuristic one only in 5 cases out of 16. As the cardinality of $|I|$ gets bigger, the heuristic procedures performs better than the warm start in most cases, and, most importantly, does not produce extremely poor results, something that the warm start is not capable of guaranteeing. For instance, in the instance of B-type patients with $|J|=3$ and $|D|=28$, the gap between the objective function of the solution of the warm start and its lower bound exceeds $100\%$. On the other hand, the gap of the heuristic w.r.t. the lower bound of the warm start has an average value of $2.7\%$.\\

\begin{table}[h!]
\centering
\scriptsize
\begin{tabular}{cccrrrrrrrrr}
\hline
\multicolumn{1}{l}{\textbf{ID}} & \multicolumn{1}{c}{\textbf{$|D|$}} & \multicolumn{1}{c}{\textbf{$|J|$}} & \multicolumn{1}{l}{\textbf{$OF$}} & \multicolumn{1}{c}{\textbf{$LB_{WS}$}} & \multicolumn{1}{c}{\textbf{$OF_{WS}$}} & \multicolumn{1}{c}{\textbf{$OF_{Heu}$}} & \multicolumn{1}{c}{\textbf{$OF_{WS}-OF$}} & \multicolumn{1}{c}{\textbf{$OF_{Heu}-OF$}} & \textbf{$\frac{OF_{WS}-LB_{WS}}{LB_{WS}}$} & \textbf{$\frac{OF_{Heu}-LB_{WS}}{LB_{WS}}$} \\
\hline
\textbf{120 A} & 14 & 2 & 19669 & 19655 & 19750 & 19892 & 81 & 223 & 0.48\% & 1.21\% \\
\textbf{120 A} & 14 & 3 & 17829 & 17767 & 17928 & 18281 & 99 & 452 & 0.91\% & 2.89\% \\
\textbf{120 A} & 28 & 2 & 22355 & 22345 & 22712 & 22547 & 357 & 192 & 1.64\% & 0.90\% \\
\textbf{120 A} & 28 & 3 & 19405 & 19148 & 21769 & 19693 & 2364 & 288 & 13.69\% & 2.85\% \\
\textbf{120 B} & 14 & 2 & 17958 & 17948 & 18097 & 18186 & 139 & 228 & 0.83\% & 1.33\% \\
\textbf{120 B} & 14 & 3 & 16516 & 16428 & 20386 & 17125 & 3870 & 609 & 24.09\% & 4.24\% \\
\textbf{120 B} & 28 & 2 & 19773 & 19693 & 21411 & 20960 & 1638 & 1187 & 8.72\% & 6.43\% \\
\textbf{120 B} & 28 & 3 & 17385 & 17242 & 40659 & 18572 & 23274 & 1187 & 135.81\% & 7.71\% \\
\textbf{120 C} & 14 & 2 & 16801 & 16792 & 16907 & 17412 & 106 & 611 & 0.68\% & 3.69\% \\
\textbf{120 C} & 14 & 3 & 15478 & 15408 & 19313 & 15838 & 3835 & 360 & 25.34\% & 2.79\% \\
\textbf{120 C} & 28 & 2 & 17734 & 17676 & 19433 & 18018 & 1699 & 284 & 9.94\% & 1.93\% \\
\textbf{120 C} & 28 & 3 & 16033 & 15865 & 16347 & 16357 & 314 & 324 & 3.04\% & 3.10\% \\
\textbf{120 D} & 14 & 2 & 19100 & 19087 & 19100 & 19100 & 0 & 0 & 0.07\% & 0.07\% \\
\textbf{120 D} & 14 & 3 & 17401 & 17334 & 17626 & 17587 & 225 & 186 & 1.68\% & 1.46\% \\
\textbf{120 D} & 28 & 2 & 21604 & 21527 & 21604 & 21604 & 0 & 0 & 0.36\% & 0.36\% \\
\textbf{120 D} & 28 & 3 & 18769 & 18613 & 20931 & 19036 & 2162 & 267 & 12.45\% & 2.27\% \\
\multicolumn{3}{c}{\textbf{Average}} & \textbf{18363} & \textbf{18283} & \textbf{20873} & \textbf{18763} & \textbf{2510} & \textbf{400} & \textbf{14.98\%} & \textbf{2.70\%} \\
\hline
\end{tabular}
\caption{Results for $|I|=120$.}
\label{tab:gaps120}
\end{table}

Optimizing the emergency (Step 7) and no-show (Step 9) back-up schedules turns out to be easy, regardless of starting from the solution obtained with the warm start solution or the heuristic. In fact, the average final gap obtained solving the problem with objective function \eqref{OF_xbar} and \eqref{OF_xhat}, respectively, is always very small (below 1\%) with the exception of a single instance where it rises to more than 50\%. It is the same instance (group D with $|I|=120$, $|J|=3$ and $|D|=28$) in which the gap of the warm start solution exceeded $100\%$.\\

In summary, the average total time it takes to generate a feasible nominal schedule and to prepare a complete set of  back-up scenarios equals to $984.3$ seconds with the heuristic procedure, and $1307.6$ seconds with the warm start. The heuristic appears to be slightly faster than the warm start. On average, starting from the heuristic solution generates better values of objective function \eqref{OF_xbar} and \eqref{OF_xhat} at the end of Step 7 and Step 9, respectively. This is due to the fact that the heuristic provides, in general, a schedule closer to the nominal one, while warm start, in the more challenging instances, may generate solution with more spare capacity. This is especially true for instances with a bigger $|I|$ cardinality.

\subsection{Case study}
\begin{figure}[h]
 \centering
\begin{subfigure}{.4\textwidth}
\centering
\includegraphics[width=.8\textwidth]{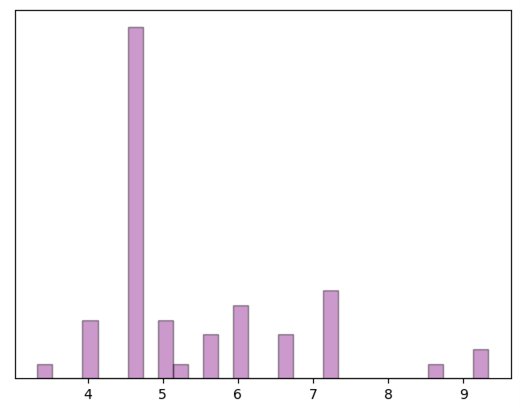} 
\caption{Distributions of surgical times in the case study analyzed, expressed as multiples of 15 minutes.}
 \label{fig:hist_buzzi}
\end{subfigure}
\hspace{3em}
\begin{subfigure}{.4\textwidth}
\centering
\includegraphics[width=.6\textwidth]{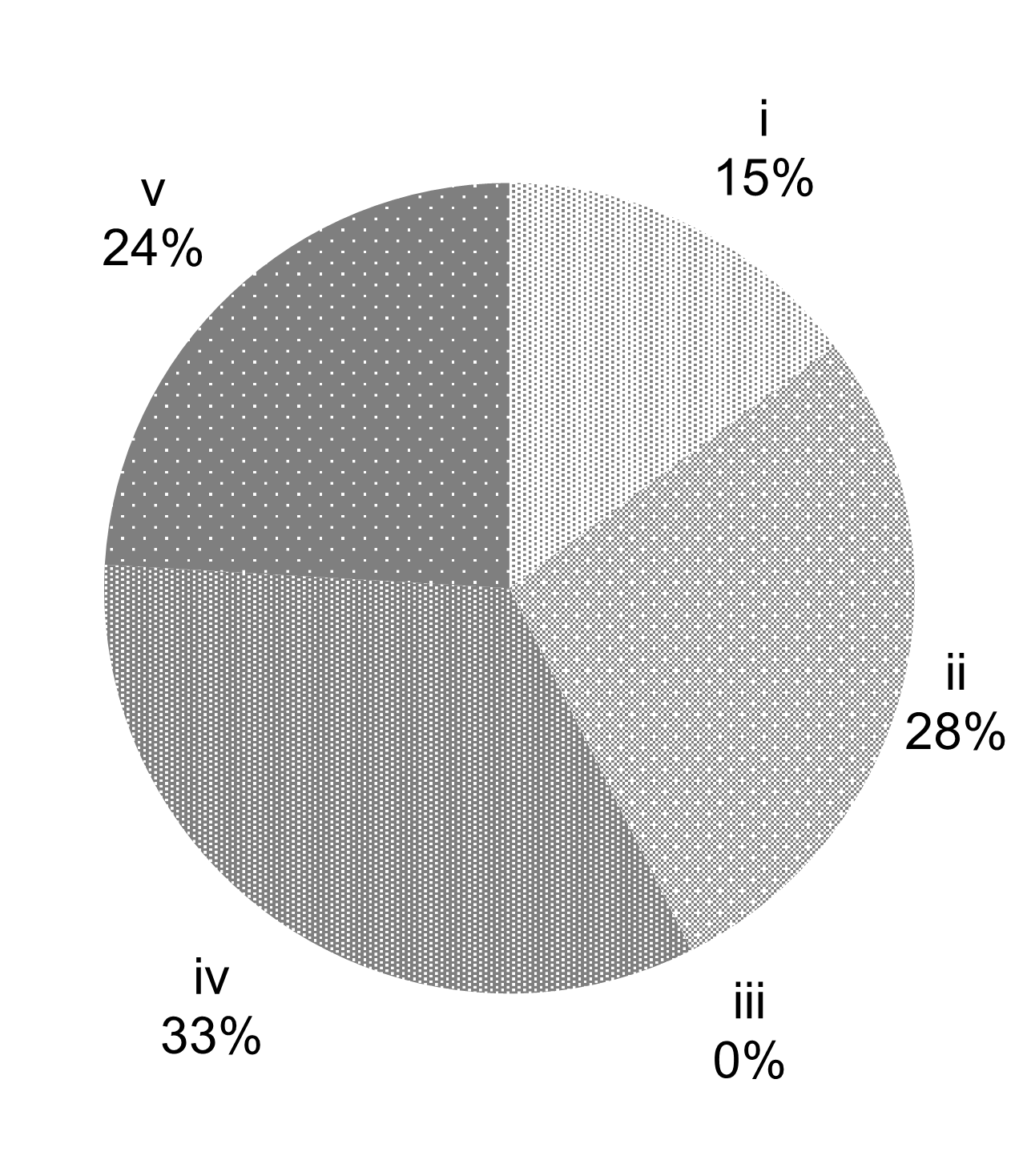} 
\caption{Distributions of the different urgency classes in the case study, \textit{i} being the most urgent.}
 \label{fig:pie_buzzi}
\end{subfigure}
\label{fig:test_buzzi}
\caption{Characterization of the case study patients}
\end{figure}

V. Buzzi Children's Hospital currently has 2 ORs, that, as described at the beginning of this manuscript, are also devoted to hosting any emergency that may arise. They are however planning to build an additional one, so, also in this case, the tested instances consider cases with $|J|=[2,3]$. Similarly, the considered time horizon is $|D|=[14,28]$, resulting in four instances. All other parameters were equivalently set as in the previous Section. This realistic instance includes $|I|=54$ real patients, whose (completely anonymous) data has been derived from historical surgeries performed in the hospital. Their surgical history has been crystallized to define a moment in which they were all in the waiting list. Their surgical times, reported in aggregate form in Figure \ref{fig:hist_buzzi}, represent therefore a realistic mix for this case study.

Table~\ref{tab:gapsbuzzi} reports about the objective functions and gaps for the Steps 1, 3 and 6. The nominal schedule problem is solved to optimality for all the four instances in an average computational time of $0.42$ seconds. The warm start always reaches the time limit, and the heuristic procedure does so in two out of four the instances, where three ORs are considered.

The results reported in Table \ref{tab:gapsbuzzi} show that, as expected, the warm start performs better than the heuristic, given that the number of patients is closer to the already tested instances with $|I| = 40$. Still, it is worth highlighting that the warm start performs way better than the heuristic. It could be due to the different distribution of urgency classes in this Case Study, which has more very urgent patient, that individually impact on the OF more. \\

\begin{table}[h]
\centering
\small
\begin{tabular}{ccrrrrrrrrr}
\hline
\multicolumn{1}{c}{\textbf{$|D|$}} & \multicolumn{1}{c}{\textbf{$|J|$}} & \multicolumn{1}{l}{\textbf{$OF$}} & \multicolumn{1}{c}{\textbf{$LB_{WS}$}} & \multicolumn{1}{c}{\textbf{$OF_{WS}$}} & \multicolumn{1}{c}{\textbf{$OF_{Heu}$}} & \multicolumn{1}{c}{\textbf{$OF_{WS}-OF$}} & \multicolumn{1}{c}{\textbf{$OF_{Heu}-OF$}} & \textbf{$\frac{OF_{WS}-LB_{WS}}{LB_{WS}}$} & \textbf{$\frac{OF_{Heu}-LB_{WS}}{LB_{WS}}$} \\
\hline
14 & 2 & 2268 & 2269 & 2278 & 2885 & 10 & 617 & 0.40\% & 27.15\% \\
14 & 3 & 2106 & 2106 & 2488 & 2367 & 382 & 261 & 18.14\% & 12.39\% \\
28 & 2 & 2400 & 2394 & 2453 & 3102 & 53 & 702 & 2.46\% & 29.57\% \\
28 & 3 & 2146 & 2138 & 2359 & 2247 & 213 & 101 & 10.34\% & 5.10\% \\
\multicolumn{2}{c}{\textbf{Average}} & \textbf{2230} & \textbf{2227} & \textbf{2395} & \textbf{2650} & \textbf{165} & \textbf{420} & \textbf{7.83\%} & \textbf{18.55\%} \\
\hline
\end{tabular}
\caption{Results for the realistic instance.}
\label{tab:gapsbuzzi}
\end{table}

In Table \ref{tab:emgynoshowbuzzi} we instead report the results associated with  back-up schedules: the warm start performs better when two ORs are considered, while heuristic provides better results when three are.\\

\begin{table}[h]
\small
\centering
\begin{tabular}{llrrrr}
\hline
\multicolumn{1}{c}{\textbf{$|D|$}} & \multicolumn{1}{c}{\textbf{$|J|$}} & \textbf{$OF^{emergency}_{WS}$} & \textbf{$OF^{emergency}_{Heu}$} & \textbf{$OF^{noshow}_{WS}$} & \textbf{$OF^{noshow}_{Heu}$} \\
\hline
14 & 2 & 2289.0 & 2891.5 & 2284.6 & 2862.7 \\
14 & 3 & 2444.0 & 2367.2 & 2433.8 & 2354.3 \\
28 & 2 & 2461.1 & 3107.2 & 2455.2 & 3105.6 \\
28 & 3 & 2359.0 & 2247.0 & 2259.1 & 2251.6 \\
\multicolumn{2}{c}{\textbf{Average}} & \textbf{2388.3} & \textbf{2653.2} & \textbf{2358.2} & \textbf{2643.5} \\
\hline
\end{tabular}
\caption{Results for emergency and no-shows scenarios, comparing the average objective function value of the solutions obtained with the warm start and with the heuristic in the realistic instances.}
\label{tab:emgynoshowbuzzi}
\end{table}

\begin{sidewaysfigure}[hptb]
\centering
\begin{subfigure}{0.49\textwidth}
 \includegraphics[width=1\textwidth]{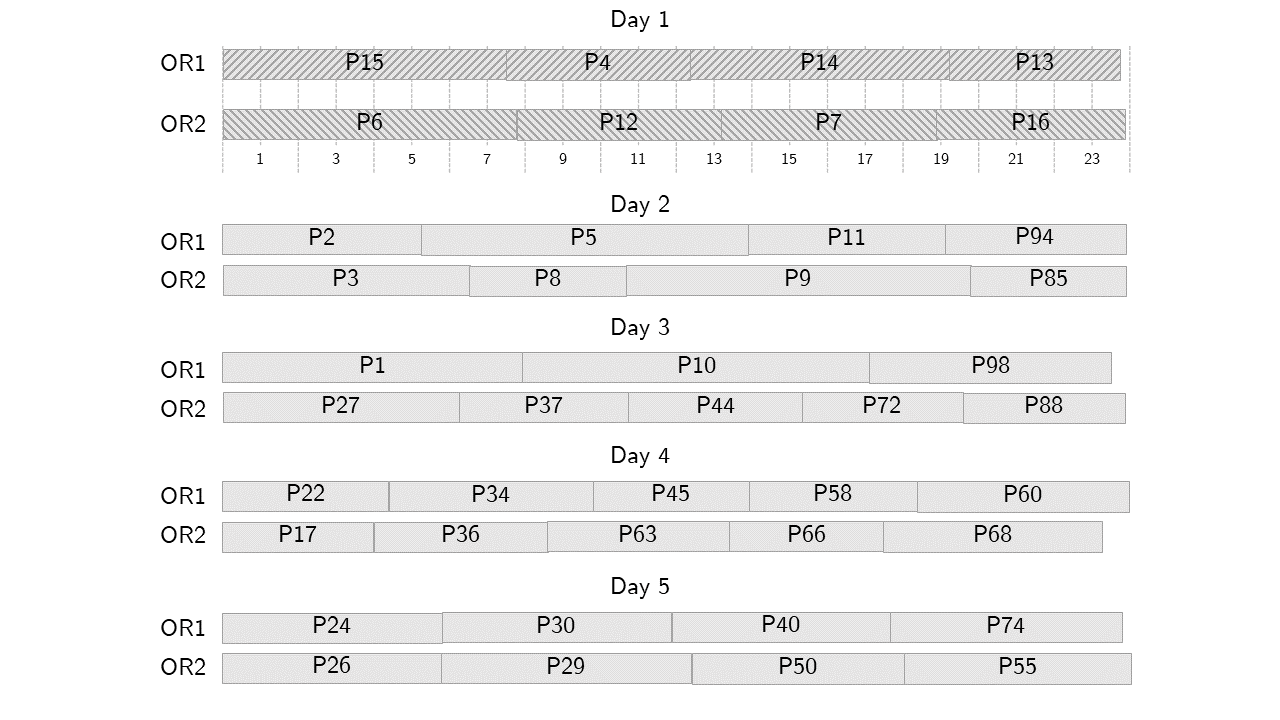}
 \caption{Sample nominal schedule}
 \label{img:nominal}
\end{subfigure}
\begin{subfigure}{0.49\textwidth}
 \includegraphics[width=1\textwidth]{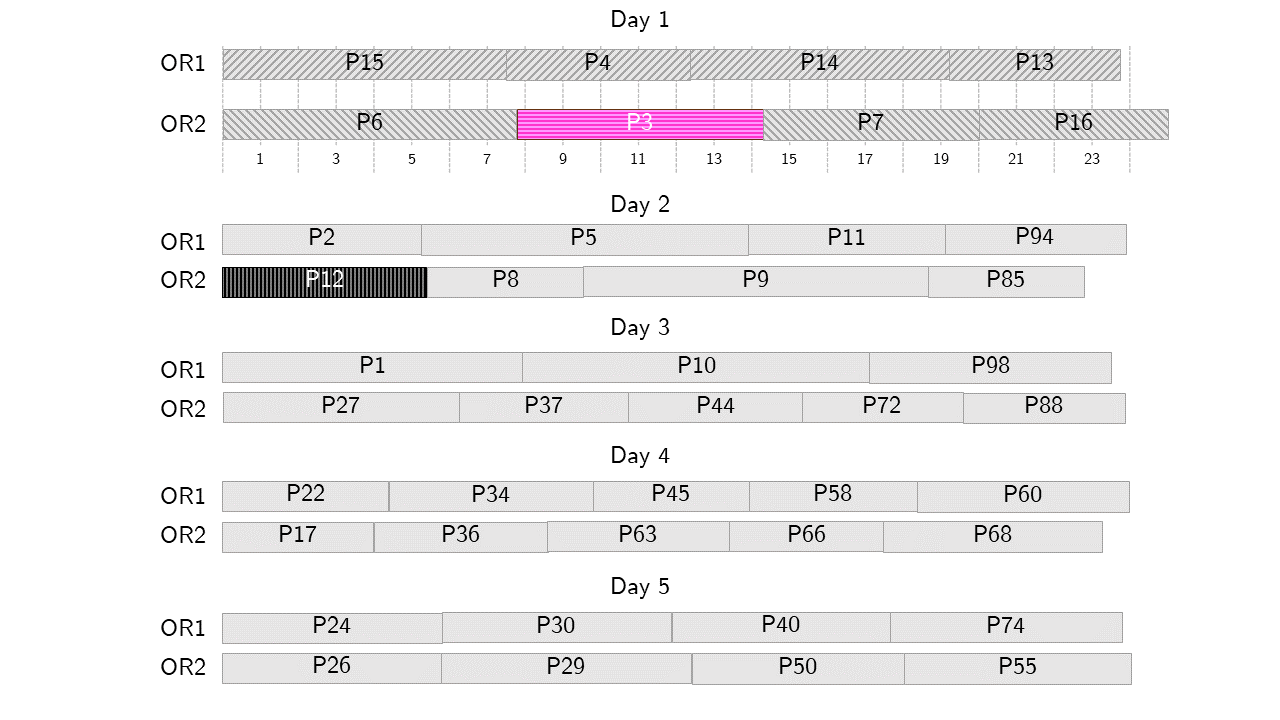}
 \caption{Schedule change for a no-show}
 \label{img:noshow}
\end{subfigure}
\begin{subfigure}{0.49\textwidth}
 \includegraphics[width=1\textwidth]{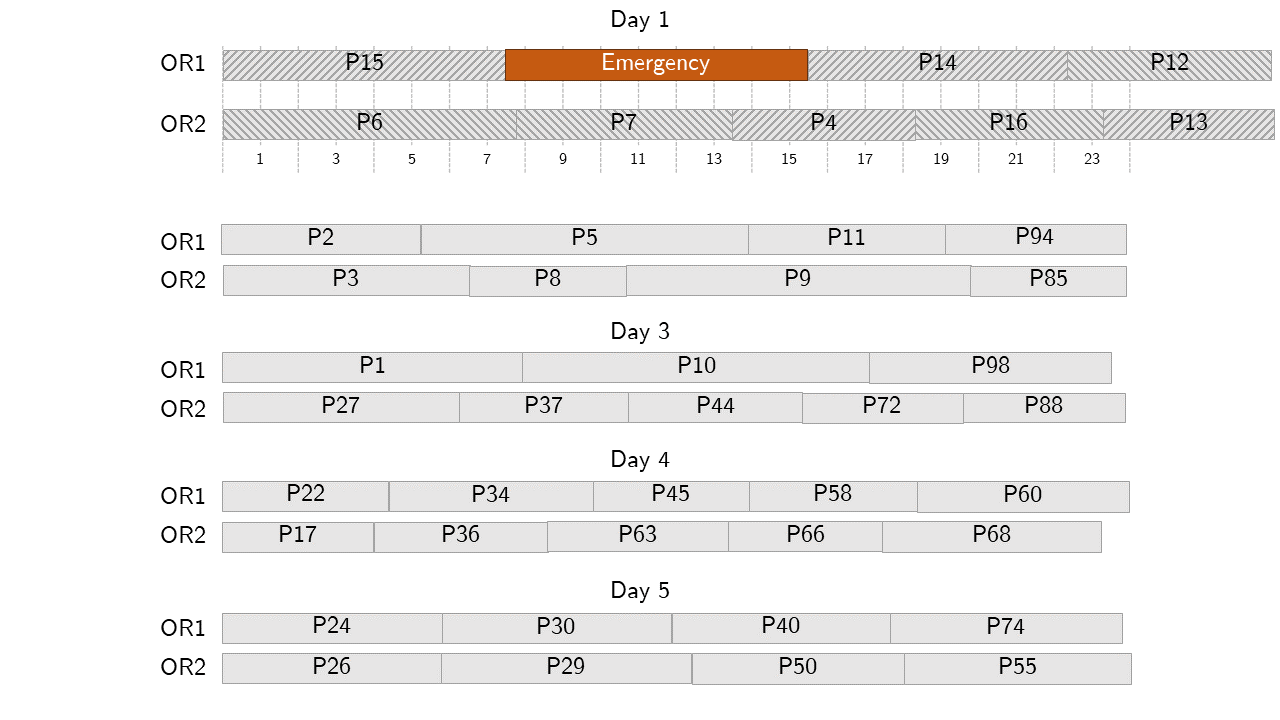}
 \caption{Schedule change in a short emergency scenario.}
 \label{img:short_emgy}
\end{subfigure}
\begin{subfigure}{0.49\textwidth}
 \includegraphics[width=1\textwidth]{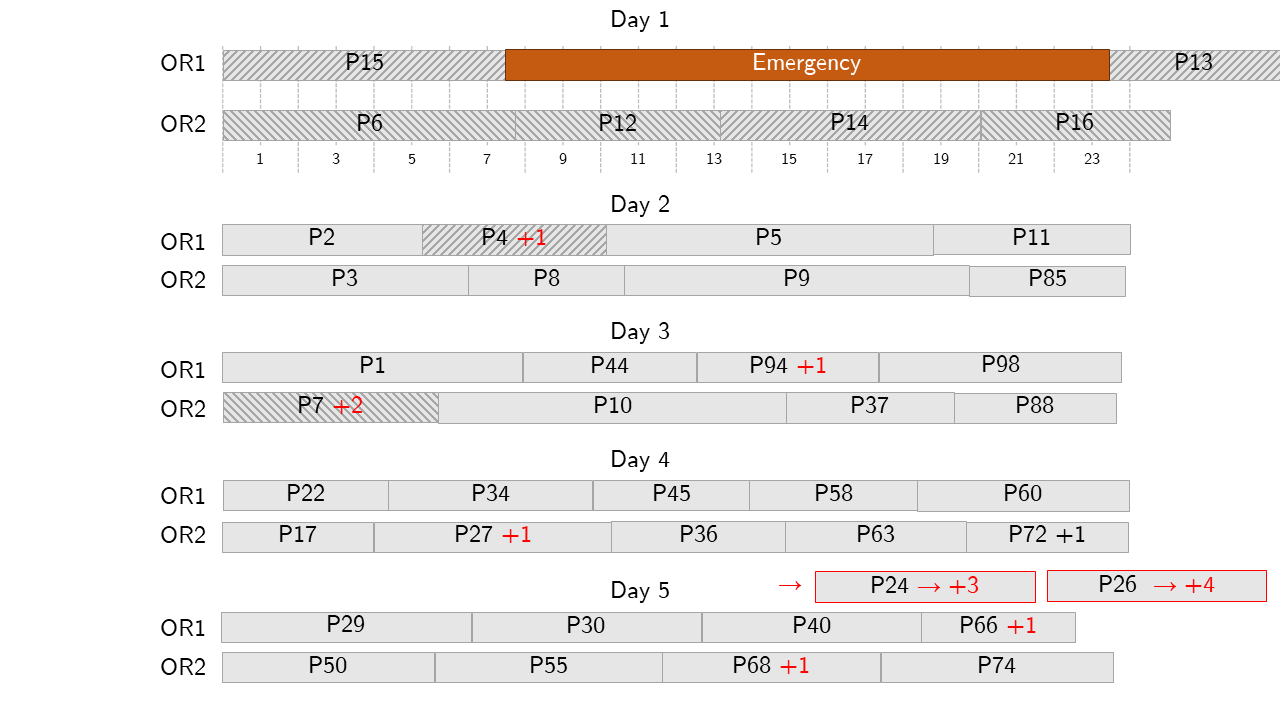}
 \caption{Schedule change in a long emergency scenario.}
 \label{img:long_emgy}
\end{subfigure}
\caption{Sample nominal schedule and two possible reschedulings.}
\label{img:example}
\end{sidewaysfigure}

Figure \ref{img:example} show graphically a set of schedules that demonstrate the functioning of the methodology. Figure \ref{img:nominal} shows an excerpt of a sample nominal schedule, showing only the first five days of the planning period. Figure \ref{img:noshow} shows what happens in case of a no-show on the first day, in this case the no-show of a patient P12 (in black/vertically striped), originally scheduled as the second patient in OR 2 in Day 1. They are substituted by back-up patient P3 (in pink/horizontally striped), who was originally scheduled on Day 2. Figures \ref{img:short_emgy} and \ref{img:long_emgy} show two different behaviours in case of emergency. In the former case, the emergency has a shorter duration, and it is possible to keep the original schedule (only moving some patients from OR1 to OR2 and vice versa) thanks to the overtime slack, parameterized by $\Omega$. In the latter case, the length of the emergency is such that the allowed overtime is not sufficient to reabsorb the disruption caused by the emergency in a single day, and, in a cascade effect, some patients are moved forward in the scheduling. The (red) number reported to the right of their ID counts how many days they have been displaced in the future. On Day 5, two patients (P24 and P26, bordered in red) exit what is showed in this excerpt, and they are rescheduled for the following period. Consider that Days 6 and 7 were modelled as non-working days.

\section{Conclusions}
We considered the problem of scheduling elective surgeries in a Children's Hospital, where disruption may occur due to emergencies and no-shows. We formulate the problem of defining a nominal schedule and a set of  back-up schedules, one for each considered disruption scenario as a disruption/restoration-based MILP model. Solving this complete model does not provide acceptable solutions even with small optimality gaps, so we implemented a warm start procedure and a heuristic procedure to provide better solutions to the problem. The warm start procedure performs better than the heuristic one in smaller instances, but becomes unstable in bigger ones. Nevertheless, both proposed approaches are able to solve the problem in reasonable time. \\
The approaches offer the great advantage of providing a set of back-up plans that are ready to be deployed immediately.\\
Future work includes embedding the approach in a rolling horizon framework, and deploying it using open source software in an actual pediatric ward.

\bibliographystyle{unsrtnat}
\bibliography{bib.bib}

\end{document}